\documentclass{amsart}

\usepackage{hyperref} 
\DeclareMathOperator{\trace}{trace} 
\DeclareMathOperator{\divergence}{div} 
\title[The Noether Theorems]{The Noether Theorems and their Application to Variational Problems on a Hyperbolic Surface} 
\newtheorem{proposition}{Proposition} 
\newtheorem{theorem}[proposition]{Theorem} 
\newtheorem{corollary}[proposition]{Corollary}

\author{Karen K.~Uhlenbeck}

\usepackage{tikz} 
\usetikzlibrary{calc,intersections,through, backgrounds,arrows,decorations.pathmorphing,backgrounds,positioning,fit,petri} 
\usetikzlibrary{calc}

\newcommand*\circled[1]{\tikz[baseline=(char.base)]{ 
\node[shape=circle,draw,inner sep=1pt] (char) {#1};}} 

\begin{document}

\maketitle

The fundamental theorem of Emmy Noether in the calculus of variations has been important historically. She came to work in G\"ottingen in 1915 at the invitation of Hilbert and Klein. In the years 1915--1918, there was a healthy competition, particularly between Einstein and Hilbert, as the outlines of general relativity took shape. General relativity is correctly attributed to Einstein, but Hilbert made substantial contributions, particularly in his development of the variational principle~\cite{H}. One topic puzzled a number of mathematicians and physicists, Hilbert in particular. Why was there no conserved quantity corresponding. to energy in general relativity? Noether was already known for her work on invariants, and was recruited by Hilbert and Klein to work on this problem. Noether surprisingly quickly produced a paper \cite{N} settling the question so thoroughly that there is even today very little to be added. At the time, there was some debate about the importance of the variational approach, and the breadth of the applicability was not at all recognized. However in the last 50 years, her theorem has morphed into a widely recognized principle embedded in the philosophy of mathematics and physics: symmetries are associated with conservation laws. But it is still useful to look at the theorem itself. An excellent and complete reference containing both the history and the mathematics can be found in \cite{K-S}.

Her original theorem is strictly a theorem in the calculus of variations and has two distinct parts which are not always differentiated. 

``Part I: If an integral $J$ in the calculus of variations is invariant under a [group] $G_\rho$ (of dimension $\rho$), there are $\rho$ linearly independent combinations among the Lagrangian expressions and their derivatives which become divergence free, and conversely\ldots .''

``Part II: If an integral $J$ in the calculus of variations is invariant under a [group] $G_{(\infty, \rho)}$ depending on $\rho$ arbitrary functions and their derivatives up to order $\sigma$, then then there are $\rho$ linearly independent identities among the Lagrangian expressions and their derivatives up to order $\sigma$, and conversely\dots .''

In part one of the theorem, each symmetry produces for every solution of the Euler-Lagrange equations, a local vector field $V$ with $\divergence V = 0$ in any Euclidean coordinates one might choose to write the integral and functions, provided the symmetry is properly expressed in those coordinates. 

Noether was nothing if not a generalist. The integrals under consideration depend on $m$ dependent functions $y^k = f^k(x)$ and their derivatives ${y^k_\alpha} ={\partial_\alpha f^k(x)}$, through order $s \geq 1$ defined on a domain $x = (x^1, \cdots ,x^n)$ in $\mathbb R^n$. The integrals are multiple integrals with an integrand, known as a Lagrangian density of the form
\[ L(x, f(x), \{\partial_\alpha f \}) (dx)^n. \]
The importance of the integrals lies in the system of Euler-Lagrange equations $1 \leq k \leq m$ for a function $u = (u^1,\cdots ,u^m)$ where the variation of $J$ vanishes. In particular, it could be a function at which $J$ takes on a minimum, but most often only the partial derivatives or variation of $J$ vanishes.
\[ \sum_\beta (-1^{|\beta|} \partial_\beta [L_{(y^k_\beta)}(x, u,\cdots, \{\partial_\alpha u\})] = 0.\]
The Lagrangian expressions are the partial derivatives of $L$ in the various directions $y^k_\beta$, but evaluated at the solution $f =u$ of the Euler-Lagrange equations.

The equations of general relativity require the number of derivatives $s = 2$. The functions on which the Einstein-Hilbert functional is defined are a metric tensor $g$ and the Lagrangian density is the scalar curvature times the volume form induced by the metric.

We are not generalists and for the purposes of this note, we will restrict to $s = 1$ and consider only one derivative. The examples in this note all have $s = 1$: hence the Lagrangian densities depend only on functions and their first derivatives, which simplifies any discussion. For Noether, a symmetry is an infinitesimal symmetry which leaves the Lagrangian density $L(x,f(x),\partial f(x))(dx)^n$ invariant (or invariant up to a divergence). It is not necessarily a symmetry preserving the functional. 

It is important to realize in the recipe that the symmetries need only be local. For a problem in the calculus of variations based on $t$ in $[0,1]$ that does not contain $t$ explicitly (like $J(f) = \int_0^1 L(f'(t),f(t)\, dt$), translation in $t$ does not preserve $[0,1]$ but it is still a local symmetry and Noether's theorem applies. How it works can be seen in her proof. Insert in the integral as the first variation the infinitesimal symmetry (the Lie algebra) acting on the solution to the Euler-Lagrange equations but multiplied by an arbitrary function $\phi$ which has local support. This vanishes because we are at a solution of the Euler-Lagrange equations. Because of the symmetry, only terms $\partial\phi$ (and higher derivatives in the case the integral depends on more than one derivative) are in the expression. There is no term with just $\phi$. In the case of many derivatives, we freely integrate by parts until an expression like the one present for only one derivative
\[ 0 = \int \sum_j V^j d_{x^j}\phi\, (dx)^n \]
for all $\phi$ with support in a ball pops out. Integrate by parts. Because $\phi$ is arbitrary, what used to be called the fundamental theorem of the calculus of variations gives us $\divergence V = 0$. This is due to the lack of a term containing only $\phi$, which in turn is due to the symmetry.

The importance of the divergence free vector fields associated with the solution of the Euler-Lagrange equations is the divergence theorem. If $\Omega$ in $\mathbb R^n$, boundary $\Omega = \mathcal S$, then $\int_{\mathcal S} [\nu,V] = \int_\Omega \divergence V$, where $\nu$ is the outward unit normal.

Our modern language interprets this in terms of forms, and since Noether's arguments are entirely local, the theorem easily extends to manifolds. The Lagrangian density is an $n$-form, the vector field $V$ is transformed into an $n-1$ form via $\theta = V\ {\lrcorner}\ (dx)^n$, and the divergence free vector fields V become closed $(n-1)$ forms $\theta$. In turn, $\theta = d\xi$ locally for an $n-2$ form. Later we will be particularly interested in $n=2$, where $\xi$ is a function.

Conservation laws for energy and momentum were well-known for many geometry and physics problems, in particular those of special relativity. A simple example is an integral in dimension $n = 1$ with the Lagrangian density $L(f, \frac{d}{dt} f)\, dt$. As mentioned before, translation in time is a local symmetry. At a solution of the Euler-Lagrange equation $u$, the quantity
\[ E(u) = \sum_k \frac{d}{dt} L_{y^k_1} \left(u,\frac{d}{dt} u\right)\frac{d}{dt} u^k - L\left(u,\frac{d}{dt} u\right) \]
is constant (in $t$). In higher dimensions, Lagrangian densities which do not depend on $x$ yield for every solution of the Euler-Lagrange equations an energy-momentum tensor $S$: an $n \times n$ symmetric matrix with divergence free rows and columns. Noether's contribution was to show that the existence of this energy-momentum tensor is due to translation invariance in time for energy and in space for momentum (to say nothing of the additional interpretation of generalized momenta coming from symmetries in the target) and fits into a larger picture of symmetry and conservation laws. It is not clear how well this was understood in the physics community before the revival of interest in her work in the 60's.

The second part of Noether's theorem answers the question on why there are no conserved quantities in general relativity. The usual sources of conservation laws such as space and time translation are all embedded in the larger group of diffeomorphisms. So according to part II of Noether's theorem, the resulting ``conservation laws'' are identities, and hold whether or not one is at a critical point of the integral. Later Schouten and Struik \cite{S-S} identified her identities as applied to general relativity as equivalent to the Bianchi identities (which were already known at the time of her article). This settles the question of conservation laws in general relativity (and for that matter, in gauge theory), although this has taken some time to be recognized. See the article by Rowe~\cite{Ro} for a more detailed discussion.

The question of gauge theory is surely an interesting one. A search of the web for ``Noether, gauge theory'' found an active discussion among physicists. No mathematicians. While some of the references are to the “Noether identities” of part~II, it is difficult to see how, as these identities hold for all fields, information about specific structures can be derived from these identities. Whatever is happening at the quantum level is a result of the representations and most likely the philosophy arising from the theorem, not the theorem itself. For those curious, the classical identity that arises from applying her theorem to the Yang-Mills functional is either $D*D*F = 0$ or $\trace \trace[F,F] = 0$, both of which are self evident even to beginners in the subject.

My current interest in the subject arose in my collaboration with George Daskalopoulos in approximating “best Lipschitz” maps between hyperbolic surfaces $f:M \to N$ studied by Thurston. To obtain approximations to the maps minimizing the Lipschitz constant of $f$, or equivalently the maximum of partial derivatives of $f$, we study integrals on variational problems which depend on a Shatten-von Neumann norm. Messing around in a fixed coordinate system on the image $N$, we found closed one forms. These were the derivatives of functions we were seeking. At some point we realized that these closed one forms in two dimensions were exactly the duals of the divergence free vector fields promised by Noether from the local symmetric space structure of $N$. A complete set of 3 closed one forms, one for each symmetry, popped out and we were on off to the races. This article will appear shortly \cite{D-U-2}.

Here is a simple example. Let $u: M \to S^1$ and $J_p(u) = \int_M |du|^p*1$, the critical points of which are called $p$-harmonic functions, be the integral under consideration. The symmetries in question are the rotations of the circle. With even a small amount of knowledge of the calculus of variations, it is not too hard to see that the Euler-Lagrange equations are
\[ d*|du|^{p-2}du = 0.\]

One does not even have to refer to Noether's recipe to see that the closed $n-1$ form corresponding to rotations of the circle is
\[ \theta = d\xi = *|du|^{p-2}du \]
where $\xi$ is a $\dim M - 2$ form. It is easy to use this without crediting Noether (as we all do), but the case of $N$ instead of $S^1$ is not so transparent. We are to make good use of these $n-1 = 1$ forms in two dimensions, but understanding the higher dimensional cases is wide open.

At some point, Daskalopoulos and I idly asked, well, what about the symmetries of the domain $M$? Now, these are not uncharted waters. For $M = \Omega \subset \mathbb R^n$, the translations are all local symmetries and Noether's theorem applies to give what physicists know as conservation of linear momentum (or energy if there is a direction representing time). Without crediting Noether, geometers for the last few decades have used these as the basis for monotonicity theorems, especially for harmonic maps and solutions of the Yang-Mills equations. They apply only approximately on curved manifolds, but that is sufficient for the estimates they are used for. Also, the proof is different from Noether's as there is often insufficient regularity for hers. A discussion of this background can be found in the paper by Bernard~\cite{B}, which also presents another application in geometry of Noether's first theorem.

We can go back to the integral, now on $\Omega$ in $\mathbb R^n$.
\[ J_p(f) = \int_\Omega |df|^p (dx)^n. \]
There are more modern proofs using the diffeomorphism directly, but we give the classical argument. The recipe for finding the conservation law corresponding to translation in the $x^j$ direction turns out in this case to be the inner product of the Euler-Lagrange equations with the partial derivative of $f = u$ in the direction $x^j$ which we write $d_{x^j} u$. We get
\[ (|du|^{p-2}\, du, d(d_{x^j}u)) = (|du|^{p-2}\, du, d_{x^j} du) = 1/p\, d_{x^j}\, |du|^p.\]

This same expression is also the same as
\[ d*(|du|^{p-2}\, du\, d_{x^j} u ) - (d*(|du|^{p-2}\, du),d_{x^j} u) = d*(|du|^{p-2}\, du, d_{x^j} u) \]
because we assume the Euler-Lagrange equations. So the divergence free (momentum) vector for the $j$-th direction is
\[ S_{(i,j)} = |du|^{p-2}( d_{x^i} u,d_{x^j} u) - \frac{1}{p}\delta_{(i,j)}|du|^p. \]

It is a bit sticky taking the dual to obtain a closed $(n-1)$ form and to keep the symmetry, and dimensions greater than 2 are more difficult.

In dimension 2, we can take the dual (star) of both indices to obtain $S^*_{(1,1)} = S_{(2,2)}$, $S^*_{(2,2)} = S^*_{(1,1)}$ and $S^*_{(2,1)} =S^*_{(1,2)} = -S_{(2,1)} = -S_{(1,2)}$. Then elementary but somewhat tedious calculations based only on the local calculus fact that in $\mathbb R^2$, $d\theta = 0$ if and only if $\theta = d\xi$, results in the following:
\begin{proposition}
  If $\Omega \subset \mathbb R^n$ and $J$ is an integral on $f:\Omega \to N$ well defined on some $L_1^p$ space which is invariant under translation, then to every solution of the Euler-Lagrange equations there is a energy momentum tensor $S$ with $d*S = 0$. If $n=2$, then $S^*$ calculated above satisfies $dS^* = 0$. Moreover any symmetric tensor $S^*$ is closed (with respect to either index) if and only if there exists $E$ such that $S^*_{(i,j)} = d_id_j E$. The function $E$ is unique up to addition by a term $ax^1 + bx^2 + c$. 
\end{proposition}

This theorem is quite general and examples besides the one we give abound. Every solution of the Euler-Lagrange from an integral invariant under translation has an energy-momentum tensor, often even when the solution has minimal regularity.

So what about when the domain is a hyperbolic surface $M$?

The calculation we did above is well-known to geometers on an arbitrary manifold, symmetric or not, and can be show to yield an equation on an energy-momentum tensor $S$ for many integrals with the domain a Riemannian manifold~$M$. Such calculations are basic steps in studying harmonic maps and Yang-Mills equations, and they come from variations of the parameterization. (In general, $S$ is a symmetric tensor with values in $T^*(M)\ \circled{s}\ T^*(M)$, and may only be in $L^1$. The equations are interpreted distributionally). However, the equation is now
\[D^*S = 0. \]
Here $D$ is the covariant derivative form of the exterior differential $d$. 

To see the connection with Noether, let $\omega$ be an infinitestimal local isometry (Killing vector field). Note that $(S,\omega)$ is now a one-form, and ${}^*(S,\omega)$ is an $n-1$ form. The following is an easy calculation.
\begin{proposition}
  If the domain is a symmetric space, $D^*S = 0$ if and only if $d^*(S,\omega) = 0$ for local Killing fields $\omega$. 
\end{proposition}

An extension of this to conformally invariant integrals and conformal Killing fields appears in a number of places (without reference to Noether) \cite{B-L}, \cite{P}.

On a surface, if we construct the symmetric tensor $S^*$ from $S$ as above, this becomes
\[ DS^* = 0.\]

The covariant derivative form of the exterior differentiation of course spoils the above calculations. On a hyperbolic surface, we have 3 linearly independent vector fields $\omega_\alpha$ corresponding to the local symmetries.

The three one-forms promised by Noether are ${}^*(S,\omega_\alpha)$. It should be noted, however, that there is a linear relation among these three forms. They are NOT linearly independent. Noether's theorem predates an understanding of Lie algebras. In another elementary but tedious calculation similar to the one in $\mathbb R^2$ in coordinates in the upper half plane we find the following:
\begin{theorem}
  On a hyperbolic surface $M$, the equation $D^*S = 0$ is valid if and only if locally $S^* = \nabla d E + R g E$ where $\nabla$ is the covariant derivative, $R$ is the constant curvature and $g$ is the metric tensor. Moreover, the kernel $K$ of the operator $\nabla d + Rg$ corresponds to functions~$E$ which are the Hamiltonians for the action of the symmetry group~$\mathrm{SO}(2,1)$ with respect to the natural symplectic form on~$M$.
\end{theorem}

Note that non-constant curvature would completely spoil the `if' part, which is easily checked once the operator is found. 

We add a brief description of this kernel $K$ in coordinates $(x^1,x^2) = (x,y)$ of the upper half space $H^2$.

Three linearly independent Killing vector fields on $H^2$ are 
\begin{align*}
  \omega_0 &= x d_x + yd_y\\
  \omega_1 &= d_y\\
  \omega_2 &= (x^2 - y^2)\, d_x + 2xy\, d_y. 
\end{align*}

The kernel $K$ of the operator $\nabla d + Rg$ is generated by the three functions 
\begin{align*}
  k_0 &= x/y\\
  k_1 &= 1/y\\
  k_2 &= (x^2 + y^2)/y. 
\end{align*}
To see that $k_j$ is the Hamiltonian for $-\omega_j$ we compute in coordinates. The computation for $k_0$ goes like this:
\[ *dk_0 = 1/y*dx-x/y^2*dy = -x/y^2\, dx - 1/y\, dy \to -xd_x - yd_y.\]

The last step identifying a one-form with a tangent vector uses the fact that the metric is conformal to the Euclidean metric with conformal factor $1/y$. The computation is similar for $k_1$ and $k_2$.

The kernel $K$ is a 3-dimensional space of functions on $\tilde{M} = H^2$ whose elements are $\sum_j a_j k_j$. If $\gamma$ is a closed curve in $M$, we lift to $\gamma: [0,1] \to H^2$. Solve $S^* = \nabla d E+ RgE$ as a solution in the cover. We can check that $E(\gamma X)$ is also a solution, so $k(X) = E(\gamma X) - E(X)$ is in~$K$. This and the homotopy invariance give us a theorem that can be used to construct an affine bundle.
\begin{corollary}
  Let M be a hyperbolic surface. To every $S^*$ in $T^*(M)\ \circled{s}\ T^*(M)$ which is closed with respect to either index, $\nu: \pi_1(M) \to K$ is well-defined on $\pi_1(M) $and satisfies
  \[ \nu(\gamma_1 \gamma_2) = \nu(\gamma_1) + \nu(\gamma_2).\]
\end{corollary}

We are specifically interested in the case where $S^*$ is a measure with support on a lamination [D-U-1]. Applications to this case as well as detailed proofs of the results outlined here can be found in a forthcoming paper. 

The references include an article on Noether's role in the relativity revolution \cite{Ro} and the references on applications of Noether's theorem \cite{B-G-G,O}. I particularly recommend the reference by Kosmann-Schwarzbach for historical background and an extensive description of related articles in both physics and mathematics \cite{K-S}.

We have no idea of either the use or real meaning of any of these observations. Or an exact statement of theorem 2 in higher dimensions. Or what applications they may have, 

Sometimes however, the pursuit of useless knowledge pays off.

\bigskip

Professor Emerita, University of Texas at Austin

Distinguished Visiting Professor, Institute for Advanced Study 
\end{document}